\documentclass[11pt]{article}
\usepackage{fullpage}
\usepackage[T1]{fontenc}
\usepackage[latin1]{inputenc}

\usepackage{amsfonts,amssymb,amsmath,xspace,amsthm}
\usepackage{url}
\usepackage{array}
\usepackage{psfrag}
\usepackage{graphicx}

\newcommand{\TsP}{$T_{\mathcal{S}}(\sigma)$\xspace}
\newtheorem{rk}{Remark}
\newtheorem{prop}{Property}
\newtheorem{conj}{Conjecture}
\newtheorem{theorem}{Theorem}
\newtheorem{lemma}{Lemma}
\newtheorem{corollary}{Corollary}

\begin{document}

%% ------------------------------------------------------------------------------

  \title{Average-case analysis of perfect sorting by reversals}

  \author{Mathilde Bouvel\footnote{CNRS, Universit\'e Paris Diderot, LIAFA, Paris, France, Supported by ANR project GAMMA BLAN07-2\_195422}\ , Cedric Chauve\footnote{Department of Mathematics, Simon Fraser University, Burnaby
    (BC), Canada}\ , Marni
    Mishna$^{\fnsymbol{footnote}}$, Dominique Rossin\addtocounter{footnote}{-1}$^{\fnsymbol{footnote}}$ }

%% ------------------------------------------------------------------------------

\maketitle

\begin{abstract}
A sequence of reversals that takes a signed permutation to the identity is perfect if at no step a common interval is broken.
Determining a parsimonious perfect sequence of reversals that sorts a signed permutation is NP-hard. Here we show that, despite 
this worst-case analysis, with probability one, sorting can be done in polynomial time. Further, we find asymptotic expressions for 
the average length and number of reversals in commuting permutations, an interesting sub-class of signed permutations.

%   Perfect sorting by reversals of a signed permutation aims at sorting
%   this permutation towards the identity permutation, by a sequence of
%   reversals that do not break any common interval. It has been shown
%   recently that this problem is NP-hard, and fixed-parameter
%   algorithms have been proposed. Here, we analyze the average-case
%   behavior of the involved parameters and show that, with probability
%   1, when a permutation is long enough, it can be sorted in polynomial
%   time. We also describe average properties of perfect sorting
%   commuting permutations, a class of signed permutations that can be
%   sorted by perfect reversals scenarios where the order of reversals
%   does not matter.
\end{abstract}

%% ------------------------------------------------------------------------------
\section{Introduction}
\label{sec:intro}

The sorting of signed permutations by reversals is a simple combinatorial problem
with a direct application in genome arrangement studies.
Different sorting scenarios provide estimates for evolutionary distance and can help explain 
the differences in gene orders between two species (see \cite{bourque-genome} for example). Initially, the 
shortest sequences (parsimonious) of reversals were sought, and polynomial time
 algorithms to find such sequences were described (\cite{hannenhalli-transforming,bergeron-inversion,tannier-advances}).
Recently, biologically motivated refinements have been considered, specifically accounting for groups of genes that 
are co-localized with the different homologous genes (genes having a single common ancestor) in the genomes of different species.
These groups are likely together in the common ancestral genome, and were not disrupted during evolution, hence, we expect them to
 appear together at every step of the evolution. In terms of our combinatorial model, a group of co-localized genes is modeled by a 
{\em common interval}, that is, a collection of sequential numbers that are not broken by any reversal move. This constraint leads us 
back to the basic algorithmic problem:

\begin{quote}
 What is the smallest number of reversals required to sort a signed permutation into the identity permutation without breaking any (subset of) 
common interval?
\end{quote}

These scenarios are called
{\em perfect}~\cite{figeac-sorting}. Because of the additional constraint, it is possible that the 
smallest perfect sorting scenario is longer that the smallest scenario.

Already it is known that this refined problem is NP-hard~\cite{figeac-sorting}. However, several authors have given sub-instances 
which can be solved in polynomial time~\cite{berard-conservation,berard-perfect,diekmann-evolution}, 
and fixed parameter tractable algorithms exist~\cite{berard-perfect,berard-more}.
For example, {\em commuting permutations} are the sub-class with the striking property that the property of a scenario 
being perfect is preserved even when the sequence of reversals is reordered. Examples of commuting scenarios arise in the study 
of mammals. All of the known sub-problems can be expressed in terms of the ``strong interval tree'' associated to a permutation, and we focus 
our attention on the structure of this tree.

Recently, several works have investigated expected properties of combinatorial objects related to genomic distance computation, 
such as the breakpoint graph ~\cite{XZS07,X08,XBS08,SLRM08}. We explore this route here, but focusing on the strong interval tree, to conduct an average case analysis of perfect sorting by reversals. First, in Section \ref{sec:prime}, we prove that for large enough $n$, with probability $1$, computing a perfect reversal sorting scenario for signed permutations can be done in time polynomial in $n$, despite the fact that this is 
NP-hard. Secondly, in Section \ref{sec:commuting}, we show that in parsimonious perfect scenarios for commuting permutations of length $n$, the average number of reversals is asymptotically $1.2n$, and the average length of a reversal is $1.02\sqrt{n}$.

%% ------------------------------------------------------------------------------
\section{Preliminaries}
%\label{sec:prelim}

%% - - - - - - - - - - - - - - - - - - - - - - - - - - - - - - - - - - - - - - - 
%\subsection{Perfect sorting by reversals}
%\label{ssec:prelim_genome}

We first summarize the combinatorial and algorithmic frameworks for
perfect sorting by reversals. For a more detailed treatment, we refer
to \cite{berard-perfect}.

% ---------------------------------------------------------------------------
\paragraph{Permutations, reversals, common intervals and perfect
  scenarios.}  

A {\em signed permutation} on $[n]$ is a permutation on the set of
integers $[n]=\{1, 2, \ldots, n\}$ in which each element has a sign,
positive or negative. Negative integers are represented by placing a
bar over them.  We denote by $Id_n$ (resp. $\overline{Id_n}$) the
identity (resp. reversed identity) permutation, $(1\ 2\ldots n)$
(resp. $(\overline{n} \ldots \overline{2}\ \overline{1})$). When the
number $n$ of elements is clear from the context, we will simply write
$Id$ or $\overline{Id}$.

An {\em interval} $I$ of a signed permutation $\sigma$ on $[n]$ is a
segment of adjacent elements of $\sigma$. The {\em content} of $I$ is
the subset of $I$ defined by the absolute values of the elements of
$I$. Given $\sigma$, an interval is defined by its content and from now,
when the context is unambiguous, we identify an interval with its
content.

The {\em reversal} of an interval of a signed permutation reverses the
order of the elements of the interval, while changing their signs. If
$\sigma$ is a permutation, we denote by $\overline{\sigma}$ the permutation
obtained by reversing the complete permutation $\sigma$. A {\em scenario}
for $\sigma$ is a sequence of reversals that transforms $\sigma$ into $Id_n$ or
$\overline{Id_n}$.  The {\em length} of such a scenario is the number
of reversals it contains. The length of a reversal is the number of elements in the interval that is reversed.

Two distinct intervals $I$ and $J$ {\em commute} if their contents
trivially intersect, that is either $I \subset J$, or $J \subset I$,
or $I \cap J = \emptyset$.  If intervals $I$ and $J$ do not commute,
they {\em overlap}.  A {\em common interval} of a permutation $\sigma$
on $[n]$ is a subset of $[n]$ that is an interval in both $\sigma$ and
the identity permutation $Id_n$. The singletons and the set $\{1, 2,
\ldots, n\}$ are always common intervals called {\em trivial common
  intervals}.

A scenario $S$ for $\sigma$ is called a {\em perfect scenario} if every
reversal of $S$ commutes with every common interval of $\sigma$. A perfect
scenario of minimal length is called a {\em parsimonious perfect
  scenario}.

A permutation $\sigma$ is said to be {\em commuting} if, there exists a perfect scenario for $\sigma$ such that 
for every pair of reversals of this scenario, the corresponding intervals commute. In such a case, this property holds for every perfect scenario for $\sigma$ \cite{berard-perfect}.

% for every common
% intervals $I$ and $J$ of $\sigma$, $I$ and $J$ commute. These permutations also 
% appear as {\em separable} permutation in the literature \cite{Iba97}.

% ---------------------------------------------------------------------------
\paragraph{The strong interval tree.}

A common interval $I$ of a permutation $\sigma$ is a {\em strong interval}
of $\sigma$ if it commutes with every other common interval of $\sigma$.  

The inclusion order of the set of strong intervals defines an $n$-leaf
tree, denoted by \TsP, whose leaves are the singletons, and whose root
is the interval containing all elements of the permutation. The strong
interval tree of $\sigma$ can be computed in linear time and space
(see \cite{bergeron-computing} for example).  We call the tree \TsP
the {\emph{strong interval tree}} of $\sigma$, and we identify a
vertex of \TsP with the strong interval it represents. In a more
combinatorial context, this tree is also called {\em substitution
  decomposition tree} \cite{AA05}. If $\sigma$ is a signed permutation, the sign of every element 
of $\sigma$ is given to the corresponding leaves in \TsP.

\begin{figure*}
  \psfrag{0}[b][b][1]{\small $\{1\}$}
  \psfrag{1}[b][b][1]{\small $\{2\}$}
  \psfrag{2}[b][b][1]{\small $\{3\}$}
  \psfrag{3}[b][b][1]{\small $\{4\}$}
  \psfrag{4}[b][b][1]{\small $\{5\}$}
  \psfrag{5}[b][b][1]{\small $\{6\}$}
  \psfrag{6}[b][b][1]{\small $\{7\}$}
  \psfrag{7}[b][b][1]{\small $\{8\}$}
  \psfrag{8}[b][b][1]{\small $\{9\}$}
  \psfrag{9}[b][b][1]{\small $\{10\}$}
  \psfrag{10}[b][b][1]{\small $\{11\}$}
  \psfrag{11}[b][b][1]{\small $\{12\}$}
  \psfrag{12}[b][b][1]{\small $\{13\}$}
  \psfrag{13}[b][b][1]{\small $\{14\}$}
  \psfrag{14}[b][b][1]{\small $\{15\}$}
  \psfrag{15}[b][b][1]{\small $\{16\}$}
  \psfrag{16}[b][b][1]{\small $\{17\}$}
  \psfrag{17}[b][b][1]{\small $\{18\}$}
  \psfrag{1, 2, 3, 4}[c][c][1]{\small $\{2,3,4,5\}$}
  \psfrag{5, 6}[c][c][1]{\small $\{6,7\}$}
  \psfrag{1, 2, 3, 4, 5, 6, 7, 8}[c][c][1]{\small $\{2,3,4,5,6,7,8,9\}$}
  \psfrag{12, 13}[c][c][1]{\small $\{13,14\}$}
  \psfrag{15, 16}[c][c][1]{\small $\{16,17\}$}
  \psfrag{9, 10, 11, 12, 13}[c][c][1]{\small $\{10, 11, 12,13,14\}$}
  \psfrag{0, 1, 2, 3, 4, 5, 6, 7, 8, 9, 10, 11, 12, 13, 14, 15, 16, 17}[c][c][1]{\small $\{1,2,3,4,5,6,7, 8, 9, 10, 11, 12,
 13, 14, 15, 16, 17,18\}$}
  \begin{center}
    \includegraphics[scale=1]{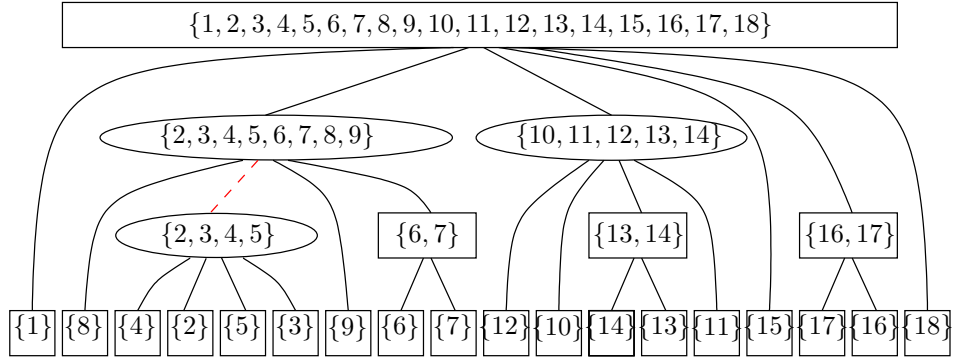}    
  \end{center}
  \caption{The strong interval tree \TsP of the permutation
    $\sigma=(1~\bar{8}~4~2~\bar{5}~3~9~\bar{6}~7~12~10~\bar{14}~13~\bar{11}~15~\bar{17}~16~18
    )$. Prime and linear vertices are distinguished by their shape.
    There are three non-trivial linear vertices, the rectangular
    vertices, and three prime vertices, the round vertices. The root
    and the vertex $\{6,7\}$ are increasing linear vertices, while the
    linear vertices $\{16,17\}$ and $\{13,14\}$ are decreasing.  }
   \label{fig:prime}
\end{figure*}

Let $I$ be a strong interval of $\sigma$ and ${\cal I}=(I_1,\dots,
I_k)$ the unique partition of the elements of $I$ into maximal strong
intervals, from left to right.  The \emph{quotient permutation} of $I$, denoted
$\sigma_I$, is defined as follows: $\sigma_I(i)$ is smaller than $\sigma_I(j)$ in $\sigma_I$ if
any element of $I_i$ is smaller (in absolute value if $\sigma$ is a signed 
permutation) than any element of $I_j$. The vertex
$I$, or equivalently the strong interval $I$ of $\sigma$, is either:
\emph{increasing linear}, if $\sigma_I$ is the identity permutation,
or \emph{decreasing linear}, if $\sigma_I$ is the reversed identity
permutation, or \emph{prime}, otherwise. For exposition purposes we
consider that an increasing vertex is positive and a decreasing vertex
is negative.  The strong interval tree as computed in the algorithm of
\cite{bergeron-computing} contains the nature -increasing/decreasing
linear or prime- of each vertex.  It can be adapted to compute also in
linear time the quotient permutation associated to each strong
interval. (See Fig.~\ref{fig:prime} for an example.)

For a vertex $I$ of $T_S(\sigma)$, we denote by $L(I)$ the set of elements of
$\sigma$ that label leaves of the subtree of $T_S(\sigma)$ rooted at $I$.

% An edge whose both vertices are prime is called a {\em prime edge}.
% The {\em prime degree} of a prime vertex $I$, denoted
% $\mathrm{deg}_p(I)$ is the number of prime edges whose two vertices
% are $I$ and a child of $I$. The {\em maximum prime degree} of \TsP is
% defined as the maximum, over all prime vertices $I$ of \TsP, of
% $\mathrm{deg}_p(I)$. 

% ---------------------------------------------------------------------------
\paragraph{The strong interval tree as a guide for perfect sorting by
  reversals.}

We describe now important properties, related to the strong interval
tree, of the algorithm described in~\cite{berard-perfect} for perfect
sorting by reversals a signed permutation. Let $\sigma$ be a signed
permutation of size $n$ and $T_S(\sigma)$
its strong interval tree, having $m$ internal vertices, called $I_1,
\ldots, I_m$, including $p$ prime vertices:

\begin{theorem}\cite{berard-perfect}\label{thm:complex}
 
\begin{enumerate}
\item The algorithm described in~\cite{berard-perfect} can compute a
  parsimonious perfect scenario for $\sigma$ in worst-case time $O(2^p
  n\sqrt{n \log(n)})$.
\item $\sigma$ is a commuting permutation if and only if $p = 0$.
\item If $\sigma$ is a commuting
  permutation, then every perfect scenario has for reversals set the 
set $\{L(I_j)|I_j \text{ has a sign different from its parent in \TsP} \}$
%  $\ell_1,\ldots,\ell_k$ be
%   the set of leaves of $T_S(\sigma')$ such that, for any $\ell_i$,
%   either $\ell_i$ is positive in $\sigma$ and its parent is decreasing
%   in $T_S(\sigma')$, or $\ell_i$ is negative in $\sigma$ and its
%   parent is increasing in $T_S(\sigma')$. Then, every parsimonious
%   perfect scenario sorting $\sigma$ has for set of reversals the set
%   $\{L(I_1),\ldots, L(I_m),\{\ell_1\},\ldots,\{\ell_k\}\}$ (reversals
%   are defined by their content here). 
\end{enumerate}

\end{theorem}

\begin{rk}
 The strong interval tree of an unsigned permutation is equivalent
to the modular decomposition tree of the corresponding labeled permutation 
graph (see \cite{berard-perfect} for example). Also commuting permutations have been investigated,
in connection with permutation patterns, under the name of {\em separable}
permutations \cite{Iba97}.
\end{rk}

% ---------------------------------------------------------------------------

%% - - - - - - - - - - - - - - - - - - - - - - - - - - - - - - - - - - - - - - - 
% \subsection{Combinatorial tools}
% \label{ssec:prelim_combi}
% 
% Description of generating functions and average properties from
% there~\cite{flajolet-analytic}.

%% ------------------------------------------------------------------------------
\section{On the number of prime vertices}
\label{sec:prime}

% The complexity of the algorithm for computing a parsimonious perfect
% scenario depends on some parameters of the strong interval tree,
% namely in our case the number of prime vertices.
Motivated by the average-time complexity of the algorithm described in \cite{berard-perfect} for
computing a parsimonious perfect scenario, we first investigate the average shape
of a strong interval tree of a permutation of size $n$. Such a tree is characterized
by the shape of the tree along with the quotient permutations labeling internal
vertices. For prime vertices, those quotient permutations correspond to {\em simple
permutations} as defined in \cite{albert-enumeration}. We first concentrate on enumerative results on
simple permutations. Next, we derive from them enumerative consequences on
the number of permutations whose strong interval tree has a given shape. Exhibiting
a family of shapes with only one prime vertex, we can prove that nearly
all permutations have a strong interval tree of this special shape.
% 
% The present section studies the average shape of a strong interval
% tree of an unsigned permutation of size $n$. Such a tree is
% characterized by the shape of the tree along with the quotient
% permutations labeling internal vertices. For prime vertices, those quotient
% permutations correspond to {\em simple} permutations as defined in
% \cite{albert-enumeration}. We first concentrate on enumerative results
% on simple permutations. Next, we derive from them enumerative
% consequences on the number of permutations whose strong interval tree
% has a given shape. Exhibiting a family of shapes with only one prime
% vertex, we can prove that nearly all permutations have a strong interval
% tree of this special shape. 
% 
% 
%% - - - - - - - - - - - - - - - - - - - - - - - - - - - - - - - - - - - - - - - 
\subsection{Combinatorial preliminaries: strong interval trees and simple permutations}\label{ssec:combPre}
Let \TsP be the strong interval tree of a permutation $\sigma$ of length $n$. From a
combinatorial point of view it is simply a plane tree (the children of a vertex are
totally ordered) with $n$ leaves and its internal vertices labeled by their quotient
permutation: an internal vertex having $k$ children can be labeled either by the
permutation $(1\ 2\ \ldots\ k)$ (increasing linear vertex), the permutation $(k\ k-1\ \ldots\ 1)$
(decreasing linear vertex) or a permutation of length $k$ whose only common
intervals are trivial (prime vertex). Due to the fact that \TsP represents the
common intervals between $\sigma$ and the identity permutation, it has two important
properties.
%\paragraph{Property 1.}
\begin{prop}
\begin{enumerate}
 \item No edge can be incident to two increasing or two decreasing
linear vertices.
\item The labeling of the leaves by the integers $\{1,\ldots , n\}$ is implicitly defined by
the labels of the internal vertices.
\end{enumerate}
\end{prop}

Permutations whose common intervals are trivial are called {\em simple permutations}.
 The shortest simple permutations are of length $4$ and are $(3\ 1\ 4\ 2)$ and
$(2\ 4\ 1\ 3)$. The enumeration of simple permutations was investigated in \cite{albert-enumeration}. The
authors prove that this enumerative sequence is not P-recursive and there is no
known closed formula for the number of simple permutations of a given size.
However, it was shown in \cite{albert-enumeration} that an asymptotic equivalent for the number $s_n$
of simple permutations of size $n$ is
\begin{equation} 
s_n = \frac{n!}{e^2} (1 -\frac{4}{n}
   +\frac{2}{n(n-1)} +\mathcal{O}(\frac{1}{n^3})) \text{ when } n
   \rightarrow \infty
\label{eqn:sn}\end{equation}
% 
% A permutation is simple if its only common intervals with the identity
% permutation are trivial. The shortest simple permutations are of
% length $4$ and are $(3\ 1\ 4\ 2)$ and $(2\ 4\ 1\ 3)$.  The enumeration
% of simple permutations is investigated in \cite{albert-enumeration}.
% The authors prove that this enumerative sequence is not P-recursive,
% so that it is very unlikely there exists some closed formula for the
% number of simple permutations of a given size. However, they obtain an
% asymptotic result:
% 
% \begin{theorem}\label{thm:asymptSimple}{\em \cite{albert-enumeration}}
%   An asymptotic equivalent for the number $s_n$ of simple permutations
%   of size $n$ is $ s_n = \frac{n!}{e^2} (1 -\frac{4}{n}
%   +\frac{2}{n(n-1)} +\mathcal{O}(\frac{1}{n^3})) \text{ when } n
%   \rightarrow \infty$
% \end{theorem}
% 
% In particular, this proves that $s_n \leq \frac{n!}{e^2}$ when $n$ is
% large enough. We refine this result by proving that:
% \begin{theorem}
%   For any $n \geq 4$, $s_n \leq \frac{n!}{e^2}$.
%   \label{thm:upper_bound_simple}
% \end{theorem}

%% - - - - - - - - - - - - - - - - - - - - - - - - - - - - - - - - - - - - - - - 
\subsection{Average shape of strong interval trees}

A \emph{twin} in a strong interval tree is a vertex of degree $2$
% labeled by $\oplus$ or $\ominus$
such that each of its two children is a leaf. A twin is then a linear
vertex. The following result, that applies
both to signed permutations and unsigned permutations, is the main result of
this section. 

\begin{theorem}\label{thm:formePermutation}
Asymptotically, with probability $1$, a random permutation $\sigma$ of size
$n$ has a strong interval tree such that the root is a prime vertex and every child
of the root is either a leaf or a twin. Moreover the probability that \TsP has
such a shape with exactly $k$ twins is $\frac{2^k}{e^2 k!}$.
%   Asymptotically, with probability $1$, a random permutation of size
%   $n$ has a strong interval tree of the form:
%   \begin{itemize}
%   \item the root is a prime vertex
%   \item every child of the root is either a leaf or a twin.
%   \end{itemize}
%   Moreover the distribution of the number $k$ of twins is given by:
%   $P(k) = \frac{2^k}{e^2 k!}$ which gives an average of $2$ twins.
\end{theorem}

% \begin{proof}
The proof follows from Lemma~\ref{lem:albert-enumeration-generalized}  and Equation \ref{eqn:sn}.
%\ref{thm:asymptSimple}.
% \qed
% \end{proof}

\begin{lemma} \label{lem:albert-enumeration-generalized}
  If $p'_{n,k}$ denotes the number of permutations of length $n$ which
  contain a common interval $I$ of length $k$ then for any fixed positive integer $c$:
  $$\sum_{k=c+2}^{n-c} \frac{p'_{n,k}}{n!} = O(n^{-c})$$
\end{lemma}

\begin{proof}
   This Lemma generalizes to any common interval the following result.

\begin{lemma}{\em \cite[Lemma 7]{albert-enumeration}}
  \label{lem:albert-enumeration} A common interval in a permutation is
  said {\em minimal} if it is not a singleton and each common interval
  included in it is trivial.  If $p_{n,k}$ denotes the number of
  permutations of length $n$ which contain a minimal common interval
  of length $k$ then for any fixed positive integer $c$:
  $$\sum_{k=c+2}^{n-c} \frac{p_{n,k}}{n!} = O(n^{-c})$$
\end{lemma}

The proof of Lemma~\ref{lem:albert-enumeration-generalized} is very similar to the article \cite{albert-enumeration}.
We have $p'_{n,k} \leq (n-k+1) k!  (n-k+1)!$. Indeed, the right hand
side counts the number of quotient permutations corresponding to $I$
($k!$), the possible values of the minimal element of $I$ ($n-k+1$)
and the structure of the rest of the permutation with one more
element which marks the insertion of $I$ ($(n-k+1)!$).  Only the
extremal terms of the sum can have magnitude ${\mathcal O}(n^{-c})$
and the remaining terms have magnitude ${\mathcal O}(n^{-c-1})$.
Since there are fewer than $n$ terms the result of Lemma
\ref{lem:albert-enumeration-generalized} follows.
%\qed
%\end{proof}

\end{proof}

\begin{proof}[Proof of Theorem~\ref{thm:formePermutation}]
  Lemma \ref{lem:albert-enumeration-generalized} with $c = 1$ gives
  that the proportion of non-simple permutations with common intervals
  of size greater than or equal to $3$ is $O(n^{-1})$. But
  permutations whose common intervals are only of size $1,2$ or $n$
  are exactly permutations whose strong interval tree has a prime root
  and every child is either a leaf or a twin.

  Then the number of permutations whose strong interval tree has a
  prime root with $k$ twins is $s_{n-k} {{n-k} \choose k} 2^k$.
 From Equation \ref{eqn:sn} the asymptotics for this number is $\frac{n! 2^k}{e^2 k!}$,
 proving
  Theorem \ref{thm:formePermutation}.
% 
%  whose
%   asymptotic for any fixed $k$ is $\frac{n! 2^k}{e^2 k!}$ proving
%   Theorem \ref{thm:formePermutation}.
  
\end{proof}
% 
% Exact experimental results for strong interval trees with up to $25$
% leaves (see Appendix) show a very fast convergence toward $1$
% (Fig.\ref{fig:exp1} in Appendix).

%% - - - - - - - - - - - - - - - - - - - - - - - - - - - - - - - - - - - - - - - 
\subsection{Average time complexity of perfect sorting by reversals}
% 
% The fact that the worst-case time complexity of the algorithm
% described in~\cite{berard-perfect} to compute a parsimonious perfect
% scenario for a signed permutation $\sigma$ of length $n$ is $O(2^p
% n\sqrt{n \log(n)})$, where $p$ is the number of prime vertices of the
% strong interval tree of the corresponding unsigned permutation,
% together with Theorem~\ref{thm:formePermutation} leads to the
% following result.

\begin{corollary}
\label{thm:complexity} The algorithm described
  in~\cite{berard-perfect} for computing a parsimonious perfect
  scenario for a random permutation runs in polynomial time with
  probability $1$ as $n \rightarrow \infty$.
\end{corollary}

\begin{proof}
Direct consequence of point 1 in Theorem~\ref{thm:complex} and of Theorem \ref{thm:formePermutation}, applied 
on signed permutations.

\end{proof}

This result however does not imply that the average complexity of this
algorithm is polynomial, as the average time complexity is the sum of
the complexity on all instances of size $n$ divided by the number of
instances. Formally, to assess the average time complexity, we need to
prove that as $n$ grows, the ratio
$$p_n = \frac{\sum_p 2^p T_{n,p}}{T_n}$$ 
is bounded by a polynomial in $n$, where $T_n$ is the number of strong
interval trees with $n$ leaves and $T_{n,p}$ the number of such trees
with $p$ prime vertices. 

\begin{figure}
  \begin{center}
    \includegraphics[scale=0.2,angle=-90]{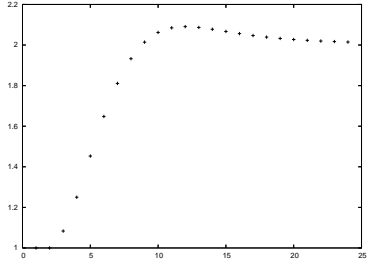}
  \end{center}
  %\vspace*{-1cm}
  \caption{$p_n$, up to $n=25$.}
  \label{fig:exp2}
\end{figure}

Let $T (x, y)$ be the bivariate generating function $T (x, y) = \sum_{k,n} T_{n,p}x^n y^p$
 Then $p_n = [x^n] T(x, 2)$. Let moreover $P(x)$ be the generating function of simple
permutations $P(x) = \sum_{n \geq 0} s_n x^n$ (whose first terms can be obtained from
entry A111111 in \cite{njas}). Using the specification for strong interval trees given in
Section \ref{ssec:combPre} and techniques described in \cite{flajolet-analytic} for example, it is immediate that
$T (x, y)$ satisfies the following system of functional equations:
$$\begin{cases}
T (x, y) = x + y  P(T (x, y)) + 2 \frac{B(x,y)^2}{1-B(x,y)}\\
B(x, y) = x + y  P(T (x, y)) + \frac{B(x, y)^2}{1 - B(x, y)}
\end{cases}
$$
By iterating these equations, we computed the $25$ first values of $p_n$ (Fig. \ref{fig:exp2})
that suggest that $p_n$ is even bounded by a constant close to $2$ and lead us to
Conjecture~\ref{conj:complexity}.

\begin{conj}\label{conj:complexity}
  The average-time complexity of the algorithm described
  in~\cite{berard-perfect} for computing a parsimonious perfect
  scenario is polynomial, bounded by ${\mathcal O}(n\sqrt n)$.
\end{conj}

%% ------------------------------------------------------------------------------
\section{Average-case properties of commuting permutations}
\label{sec:commuting}

We now study the family of commuting (signed) permutations and more precisely
the average number of reversals in a parsimonious perfect scenario for a
commuting permutation and the average length of a reversal of such a scenario.

Let $\sigma$ be a commuting permutation of size $n$, i.e. a signed permutation whose
strong interval tree \TsP has no prime vertex. It follows from the combinatorial
specification of strong interval trees given in Section  \ref{ssec:combPre} that \TsP is simply a
plane tree with internal vertices having at least two children  and a sign on the root (that  defines implicitly the signs of the
other internal vertices from point $1$ in Property $1$ and the labels $\{1 \ldots n\}$ of the leaves). These trees are then Schr\"oder
trees (entry A001003 in the On-Line Encyclopedia of Integer Sequences \cite{njas}) with
a sign on the root.

\begin{theorem}\label{thm:distance}
  The average length of a parsimonious perfect scenario for a
  commuting permutation of length $n$ is asymptotically
  $$\frac{1+\sqrt{2}}{2}n  \simeq 1.2n.
  $$
\end{theorem}

\begin{proof}
From the previous section and points 2 and 3 in Theorem \ref{thm:complex}, the problem
of computing the expected number of reversals of a parsimonious perfect scenario
reduces to computing the expected number of internal vertices of \TsP
other than the root (because two adjacent linear vertices cannot have the same
sign) and the expected number of leaves whose sign in $\sigma$ differs from the sign of
its parent in \TsP.

The expected number of leaves whose sign in $\sigma$ is different from
its parent in \TsP is obviously $n/2$, as the sign of the
leaf and of its parent are independent.

To compute the average number of internal vertices in a Schr\"oder
tree, we use symbolic methods as defined in \cite{flajolet-analytic}.
Let us define the bivariate generating function $ S(x,y) = \sum_{k,n}
S_{n,k} x^n y^k $ where $S_{n,k}$ denotes the number of Schr\"oder
trees with $n$ leaves and $k$ internal vertices. The average number of
internal vertices in a Schr\"oder tree with $n$ leaves is
$$\frac{\sum_{k} k S_{n,k}}{ \sum_{k} S_{n,k}} =
\frac{[x^n] \frac{\partial S(x,y)}{\partial y}|_{y=1}}{[x^n] S(x,1)}.
$$

A Schr\"oder tree can be recursively described as a
single leaf, or a root having at least two children, which are again
Schr\"oder trees.  Consequently, $S(x,y)$ satisfies the equation 
$$
S(x,y) = x + y \frac{S(x,y)^2}{1-S(x,y)},
$$
and solving this equation gives
\begin{eqnarray}
  S(x,y) = \frac{(x+1) - \sqrt{(x+1)^2 -4x(y+1)}}{2(y+1)}. \label{eqn:S(x,y)}
\end{eqnarray}

We compute an asymptotic equivalent of the number $[x^n] S(x,1)$, the number of 
Schr\"oder trees (\cite[entry
A001003]{njas}).

\paragraph*{Asymptotic study of $S(x,1)$.}

By Equation \ref{eqn:S(x,y)} we obtain
$$
S(x,1) = \frac{(x+1) - \sqrt{(x+1)^2 -8x}}{4} = \frac{(x+1) -
  \sqrt{(1-\frac{x}{3+2\sqrt{2}})(1-\frac{x}{3-2\sqrt{2}})}}{4},
$$

which yields the equivalent when $x \rightarrow 3-2\sqrt{2}$, $x <
3-2\sqrt{2}$
$$  S(x,1) 
% & \sim & \frac{4-2\sqrt{2}}{4} - \frac{1}{4} \sqrt{1- \frac{3-2\sqrt{2}}{3+2\sqrt{2}}} (1-\frac{x}{3-2\sqrt{2}})^{1/2}\\
   \sim  \frac{2-\sqrt{2}}{2} - \frac{\sqrt{3\sqrt{2}-4}}{2} (1-\frac{x}{3-2\sqrt{2}})^{1/2} .
$$

Applying the techniques of \cite[chapters $4$ and $6$]{flajolet-analytic}  gives the following equivalent of the coefficients $[x^n]
S(x,1)$ when $n \rightarrow \infty$:
$$
  [x^n] S(x,1) 
%& \sim & - \frac{\sqrt{3\sqrt{2}-4}}{2} (\frac{1}{3-2\sqrt{2}})^n (-\frac{1}{2\sqrt{\pi n^3}}) \\
   \sim  \frac{\sqrt{3\sqrt{2}-4}}{4} (3+2\sqrt{2})^n \frac{1}{\sqrt{\pi n^3}}.
$$

\paragraph*{Asymptotic study of $\frac{\partial S(x,y)}{\partial y}|_{y=1}$.}
By Equation \ref{eqn:S(x,y)} we obtain
$$
  \frac{\partial S(x,y)}{\partial y}|_{y=1} 
   =  \frac{(x-1)^2 - (x+1)\sqrt{(x+1)^2 -8x}}{8 \sqrt{(1-\frac{x}{3+2\sqrt{2}})(1-\frac{x}{3-2\sqrt{2}})}}.
$$

From the above expression, we can obtain an equivalent of
$\frac{\partial S(x,y)}{\partial y}|_{y=1}$ when $x \rightarrow
3-2\sqrt{2}$, $x < 3-2\sqrt{2}$. Namely,
$$
  \frac{\partial S(x,y)}{\partial y}|_{y=1} 
   \sim  \frac{3-2\sqrt{2}}{4\sqrt{3\sqrt{2}-4}} (1-\frac{x}{3-2\sqrt{2}})^{-1/2}.
$$
As before, we deduce that an equivalent of the coefficients $[x^n]
\frac{\partial S(x,y)}{\partial y}|_{y=1}$ when $n \rightarrow
\infty$ is

$$  [x^n] \frac{\partial S(x,y)}{\partial y}|_{y=1} 
   \sim  \frac{3-2\sqrt{2}}{4\sqrt{3\sqrt{2}-4}} (3+2\sqrt{2})^n \frac{1}{\sqrt{\pi n}}
$$

An equivalent of the average number of internal vertices in a
Schr\"oder tree with $n$ leaves is now easily derived as
$$
\frac{[x^n] \frac{\partial S(x,y)}{\partial y}|_{y=1}}{[x^n] S(x,1)}
\sim \frac{3-2\sqrt{2}}{3\sqrt{2}-4} n \sim \frac{n}{\sqrt{2}}.
$$

\paragraph*{Combining all results together}
The number above is the the average number of
internal vertices in Schr\"oder trees with $n$ leaves, including the root if it is not a
leaf (i.e. $n > 1$). A given Schr\"oder tree  with $n$ leaves can have
its internal vertices and leaves signed in $2^{n+1}$ ways ($2$ choices for the sign of the
root, that define the signs of all other internal vertices, and $2^n$ choices for the
signs of the $n$ leaves). As these signs do not change the number of internal vertices
of the tree, the average number of internal vertices in such signed Schr\"oder trees
does not change. We also have to discard the root as it does not define a reversal,
but this does not change the asymptotic behaviour and adding $n/2$ to account for
signed leaves that define reversals, we obtain
$$
\frac{1+\sqrt{2}}{2}n
$$

\end{proof}

\begin{rk} It is interesting to note the large representation of reversals of length
$1$, that composes almost half of the expected reversals. A similar property was
observed in \cite{lefebvre} on datasets of bacterial genomes.
\end{rk}

\begin{theorem}\label{thm:length}
  The average length of a reversal in a parsimonious perfect scenario for a
  commuting permutation of length $n$ is asymptotically 
  $$
\frac{2^{7/4} \sqrt{3-2 \sqrt 2}}{1+\sqrt{2}} \sqrt{ \pi n} \simeq 1.02 \sqrt n
  $$
\end{theorem}
\begin{proof}
We want to compute the ratio between the average  sum of
the lengths of the reversals of a parsimonious perfect scenario for a
commuting permutation and the average length of such a scenario. The
later was obtained above (Theorem~\ref{thm:distance}), and we
concentrate on the former.

A reversal defined by a vertex $x$ of the strong interval tree
\TsP is of length $L(x)$ (it reverses the segment of the
signed permutation that contains the leaves of the subtree rooted at
$x$, see~\cite{berard-perfect}). 
We first focus on the average value of the sum of the sizes of all
subtrees in a Schr\"oder tree. For simplicity in the computation, we
will also count the whole tree and the leaves as subtrees (obviously
of size $1$), which will give the same quantity we want to compute, up
to subtracting $3/2 \cdot n$ to the final result. We first define the
bivariate generating function (that we call again $S$, but which is
slightly different) following the standard analytic method defined in \cite{flajolet-analytic}
$$
S(x,y) = \sum_{k,n} S_{n,k} x^n y^k
$$
where $S_{n,k}$ denotes the number of Schr\"oder trees with $n$ leaves
and sizes of subtrees (including leaves and the whole tree) that sum
to $k$. The average value of the sum of the sizes of every subtree in
a Schr\"oder tree with $n$ leaves is
$$ \frac{\sum_{k}
  k S_{n,k}}{ \sum_{k} S_{n,k}} = \frac{[x^n] \frac{\partial
    S(x,y)}{\partial y}|_{y=1}}{[x^n] S(x,1)}.
$$

A Schr\"oder tree can be recursively described as a single leaf or a
root having at least two children, which are again Schr\"oder trees.
In the second case, the subtrees are those involved in the children of
the root, plus the tree itself (which is a subtree of size $n$), which
gives the functional equation~\ref{eqn:S(xy,y)}:
\begin{equation}
  S(x,y) = xy + \frac{S(xy,y)^2}{1-S(xy,y)}.\label{eqn:S(xy,y)}
\end{equation}

Since this equation involves both $S(x,y)$ and $S(xy,y)$, we cannot
extract from it an expression for $S(x,y)$ as in the proof 
of Theorem~\ref{thm:distance}. But since the average value of the sum of the sizes of
every subtree in a Schr\"oder tree with $n$ leaves can be obtained by
$ \frac{\sum_{k} k S_{n,k}}{ \sum_{k} S_{n,k}} = \frac{[x^n]
  \frac{\partial S(x,y)}{\partial y}|_{y=1}}{[x^n] S(x,1)}, $ we do no
need to compute $S(x,y)$ but only $S(x,1)$ and $\frac{\partial
  S(x,y)}{\partial y}|_{y=1}$.

\paragraph*{Asymptotic study of $S(x,1)$.}

By Equation \ref{eqn:S(xy,y)} we obtain $ S(x,1) = \frac{(x+1) -
  \sqrt{(x+1)^2 -8x}}{4} , $ which is the same function as in the
proof of Theorem~\ref{thm:distance}.

Hence, 
$$ [x^n] S(x,1) \sim \frac{\sqrt{3\sqrt{2}-4}}{4}
(3+2\sqrt{2})^n \frac{1}{\sqrt{\pi n^3}}.  $$

\paragraph*{Asymptotic study of $\frac{\partial S(x,y)}{\partial
    y}|_{y=1}$.}

Deriving Equation \ref{eqn:S(xy,y)} and setting $y = 1$ gives:
$$
\begin{cases}
  \frac{\partial S}{\partial x} (x,1) = 1 +  \frac{\partial S}{\partial x} (x,1) \cdot \frac{2S(x,1) -S(x,1)^2}{(1-S(x,1))^2}\\
  \frac{\partial S}{\partial y} (x,1) = x + \Big( x \frac{\partial
    S}{\partial x} (x,1) + \frac{\partial S}{\partial y} (x,1) \Big)
  \cdot \frac{2S(x,1) -S(x,1)^2}{(1-S(x,1))^2}.
\end{cases}
$$
From this system, we can extract the following equation where $S(x,1)$ has been computed before:
$$
\frac{\partial S(x,y)}{\partial y}|_{y=1} = \frac{\partial S}{\partial
  y} (x,1) = \frac{x}{(1-C)^2}, \textrm{ where } C = \frac{2S(x,1)-S(x,1)^2}{(1-S(x,1))^2}.
$$

The singularity closest to the origin is $3-2\sqrt{2}$, and the Taylor
development of the above around this singularity gives:
$$
\frac{\partial S(x,y)}{\partial y}|_{y=1} \sim \frac{3-2\sqrt 2}{2(1-\frac{x}{3 - 2\sqrt 2})}
$$

Applying the techniques of \cite{flajolet-analytic}, this yields the
following equivalent of the coefficients $[x^n] \frac{\partial
  S(x,y)}{\partial y}|_{y=1}$ when $n \rightarrow \infty$:
$$  [x^n] \frac{\partial S(x,y)}{\partial y}|_{y=1} 
\sim \frac{3-2\sqrt 2}{2}(3+2\sqrt{2})^n 
$$

Then 
$$
\frac{[x^n] \frac{\partial S(x,y)}{\partial y}|_{y=1}}{[x^n] S(x,1)}
\sim 2^{3/4} \sqrt{3-2\sqrt 2}\ \sqrt{\pi n^3}.
$$
gives the average sum of the sizes of all subtrees of a Schr\"oder
tree.

This is independent of the signs added to give the strong interval
tree of a commuting permutation, so this number is also the expected
sum of the sizes of all subtrees of a the strong interval tree
associated to a random commuting permutation. To get the expected sum
of the lengths of the reversals of a parsimonious perfect scenario for a random
commuting permutation, we need to remove the size of the whole tree,
that was counted as a subtree ($n$), the size of the $n$ subtrees
defined by the leaves ($n$) and to add the contribution of the
reversals of size $1$ ($n/2$ on the average), which does not change the above asymptotics.

% 
% 
% An equivalent of the average value of the sum of the sizes of all
% subtrees in a Schr\"oder tree with $n$ leaves is now easily derived as
% $$
% \frac{[x^n] \frac{\partial S(x,y)}{\partial y}|_{y=1}}{[x^n] S(x,1)}
% \sim 2^{3/4} \sqrt{3-2\sqrt 2}\ \sqrt{\pi n^3}.
% $$
% 
% The above result does take into account the whole tree and all leaves, that should not be 
% counted but these terms are negligible asymptotically. Hence, the average sum of the lengths of the reversals of a
% parsimonious perfect scenario for a commuting permutation of size $n$
% is asymptotically 
% $$
% 2^{3/4} \sqrt{3-2\sqrt 2}\ \sqrt{\pi n^3}.
% $$

Dividing by the average number of reversals of such a scenario
(Theorem~\ref{thm:distance}), we obtain Theorem~\ref{thm:length}.

%%%%%%%%%%%%%%%%%%%%%%%%%%%%%%%%%%%%%%%%%%%%%%%%%%%%%%%%%%%%%%%%%%%%%%%%%
%%%%%%%%%%%%%%%%%%%%%%%%%%%%%%%%%%TODO 
%%%%%%%%%%%%%%%%%%%%%%%%%%%%%%%%%%%%%%%%%%%%%%%%%%%%%%%%%%%%%%%%%%%%%

\end{proof}

%% ------------------------------------------------------------------------------
\section{Conclusion}
\label{sec:conclusion}

We showed that perfect sorting by reversals, although an intractable
problem, is very likely to be solved in polynomial time for random
signed permutations. This result relies on a study of the shape of a
random strong interval tree that shows that asymptotically such trees
are mostly composed of a large prime vertex at the root and small
subtrees. As the strong interval tree of a permutation is equivalent
to the modular decomposition tree of the corresponding labeled
permutation graph~\cite{berard-perfect}, this result agrees with the
general belief that the modular decomposition tree of a random graph
 has a large prime root. We were also able to give precise asymptotic 
results for the expected lengths of a parsimonious perfect scenario and
of a reversal of such a scenario for random commuting permutations. 

Our research leaves at least one open problem: proving that computing
a parsimonious perfect scenario can be done in polynomial time on the
average. It would also be interesting to see if our approach can be
extended to the perfect rearrangement problem for the
Double-Cut-and-Join model that has been introduced
recently~\cite{berard-perfect2} and has the intriguing property that
instances that were hard to solve for reversals are can be solved
in polynomial time in the DCJ context and conversely.

%% ------------------------------------------------------------------------------
\bibliographystyle{abbrv} \bibliography{biblio-short}

%% ------------------------------------------------------------------------------

\end{document}